\documentclass[12pt,reqno]{amsart}
\usepackage{mathrsfs}
\usepackage{amsgen}
\usepackage{amscd}
\usepackage{amsmath}
\usepackage{latexsym}
\usepackage{amsfonts}
\usepackage{amssymb}
\usepackage{amsthm}
\usepackage{graphicx}
\usepackage{color}
\usepackage{amsthm,amssymb}
\usepackage{graphicx}
\usepackage{enumerate}
\usepackage{amssymb,amsmath,graphicx,amsfonts,euscript}
\usepackage{color}
\usepackage{mathrsfs}
\usepackage{empheq}
\usepackage{cases}
\usepackage{cite}
\usepackage{bbm}
\usepackage[colorlinks=true,backref]{hyperref}

\usepackage[margin=3cm, a4paper]{geometry}

\def\H{{\cal H}}

\def\N{\mathbb{N}}
\def\R{\mathbb{R}}
\def\Z{\mathbb{Z}}
\def\T{\mathbb{T}}

\def\H2{H^2(\R^N)}
\def\L2{L^2(\R^N)}

\def\norm#1{\left\|#1\right\|}

\def\normb#1{\big\|#1\big\|}

\def\fbrk#1{\left\lbrace#1\right\rbrace}

\def\jb#1{\langle#1\rangle}

\newcommand{\hra}{{\hookrightarrow}}

\def\loe{\leqslant}
\def\goe{\geqslant}
\def\lsm{\lesssim}

\def\ep{\varepsilon}

\def\H1{H^1(\R)}

\newcommand{\al}{\alpha} 
\newcommand{\ga}{\gamma} 
 
\newcommand{\om}{\omega} \newcommand{\la}{\lambda}

\newcommand{\De}{\Delta} \newcommand{\Om}{\Omega}

   \newcommand{\I}{\infty}

\newcommand{\EQ}[1]{\begin{align*}\begin{split} #1 \end{split}\end{align*}}
\newcommand{\EQn}[1]{\begin{align}\begin{split} #1 \end{split}\end{align}}

 \newcommand{\Del}[1]{}

\numberwithin{equation}{section}

\newtheorem{thm}{Theorem}[section]

\newtheorem{lem}[thm]{Lemma}

\newtheorem{definition}[thm]{Definition}
\theoremstyle{remark}

\newtheorem*{exam*}{Examples}

\allowdisplaybreaks[4]

\begin{document}

	\setcounter{page}{1}

	\title[Incompressible Navier-Stokes equations]{Remarks on the almost sure well-posedness for incompressible Navier-Stokes equations with arbitrary regularity}
	
	\author{Ruobing Bai}
	\address{Center for Applied Mathematics\\
		Tianjin University\\
		Tianjin 300072, China}
	\email{baimaths@hotmail.com}
	\thanks{}
	
	\author{Jia Shen}
	\address{Center for Applied Mathematics\\
		Tianjin University\\
		Tianjin 300072, China}
	\email{shenjia@tju.edu.cn}
	\thanks{}

	
	

	\begin{abstract}\noindent
		In this note, we study the random data problem for incompressible Navier-Stokes equations in Euclidean space with arbitrary regularity, using the randomization introduced in \cite{Shen-Soffer-Wu-2021-1}. 
	\end{abstract}
	
	\maketitle

	\vskip 0.2cm
	\section{Introduction}
	In this note, we consider the Cauchy problem for incompressible Navier-Stokes equations
	\begin{equation}\label{NS}
		\left\{ \aligned
		&\partial_t u-\Delta u +(u\cdot \nabla) u+\nabla p=0,
		\\
		&\mbox{div } u=0,
		\\
		&u(0,x)=u_0(x),
		\endaligned
		\right.
	\end{equation}
	where $u(t,x):\R^+\times \R^d\rightarrow \R^d$, $p(t,x):\R^+\times \R^d\rightarrow \R$, and $d\goe 2$. This system describes the flow of viscous incompressible fluid, where $u$ and $p$ denote the fluid velocity and the pressure, respectively.
	
	The Navier-Stokes equations \eqref{NS} is invariant under the scaling transform: for any $\la>0$,
	\EQ{
		u_{\lambda}(t,x)=\lambda u(\lambda^2t,\lambda x)\text{, and }p_{\lambda}(t,x)=\lambda^2 p(\lambda^2t,\lambda x).
	}
	Then, this symmetry gives the critical spaces
	\EQ{
		\dot H_x^{\frac d2-1}(\R^d)\hookrightarrow L_x^d(\R^d) \hookrightarrow \dot B_{p,q}^{\frac dp-1}(\R^d) \hookrightarrow \dot B_{p,\I}^{\frac dp-1}(\R^d) \hookrightarrow BMO^{-1} \hookrightarrow \dot B_{\I,\I}^{-1},
	}
	for any $2\loe p< \infty$ and $2\loe q<\infty$.

	Now, we briefly review the mathematical theory for the Navier-Stokes equations. The theoretical study originates from the pioneering works of Leray \cite{Leray-1934} and Hopf \cite{Hopf-1951} about the existence of global weak solution. Since then, there has been much work on the uniqueness and regularity of weak solutions. The partial regularity for a class of Leray-Hopf solutions to the Navier-Stokes equations was developed by Scheffer \cite{Scheffer-1977}, who introduced the local energy inequality. Then, Caffarelli, Kohn and Nirenberg \cite{Caffarelli-Kohn-Nirenberg-1982} proved the existence of suitable weak solutions and that 1-dimensional Hausdorff measure of the singular set is zero. More results for partial regularity theory can be referred to \cite{Escauriaza-Seregin-Sverak-2003, Lin-1998, Koch-Nadirashvili-Seregin-Sverak-2009} and the references therein. Recently, Buckmaster and Vicol \cite{Buckmaster-Vicol-2019} proved the nonuniqueness of weak solution in the class of weak solutions with finite kinetic energy via the convex integration method.
	
	Next, we recall the mild solution results for the Navier-Stokes equations. Fujita and Kato \cite{Fujita-Kato-1964} first obtained the local well-posedness in ${H}^{s}(\R^3)$ when $s\goe \frac 1 2$ and small data global well-posedness in ${H}^{\frac 1 2}(\R^3)$. In \cite{Chemin-1992}, Chemin proved the small data global well-posedness in ${H}^{s}(\R^3)$ when $s>\frac 1 2$. The related results for $L^p(\R^d)$ were given by Kato \cite{Kato-1984}. More critical spaces for the global existence of mild solutions with small initial data have been considered, such as the space  $\dot{B}_{p, \infty}^{d/p -1}$ by Cannone \cite{Cannone-1997}, the space $BMO^{-1}$ by Koch-Tataru \cite{Koch-Tataru-2001} and the space $\chi^{-1}$ by Lei-Lin \cite{Lei-Lin-2011}. However, Bourgain and Pavlovi\'{c} \cite{Bourgain-Pavlovic-2008} obtained the ill-posedness of Navier stokes equations in the largest critical space $\dot{B}_{\infty,\infty}^{-1}(\R^3)$. Moreover, even though the spaces $\dot B^{-1}_{\I,q}$ is embedded into $BMO^{-1}$, Wang \cite{Wang15Adv} proved that the equation is ill-posed in such space. Some other results can be found in \cite{Chemin-Gallagher-2009, Lemarie-Rieusset-2002} and the references therein.
	
	The existence of another class of solutions called self-similar solutions has also been studied. 
	The null solution is a backward self-similar solution of the Navier-Stokes equations, and for other cases see the discussions in \cite{Leray-1934,Tsa98ARMA,NRS96Acta,Escauriaza-Seregin-Sverak-2003}. The forward self-similar solutions were constructed in small data setting \cite{Cannone-Planchon-1996, Chemin-1999, Koch-Tataru-2001}, and for the large data case, see Jia-\v{S}ver\'{a}k's result \cite{Jia-Sverak-2014}.
	
	While currently the global regularity of 3D Navier-Stokes equation without any size restriction is unsloved, the related 2D result has been obtained by Ladyzhenskaya \cite{Ladyzhenskaya-1963}. In fact, in the 2D case, the critical space is $L_x^2$ and the energy inequality provides an a priori control at the same scaling. The global result was extended by Gallagher and Planchon to more general critical Besov spaces, see \cite{GP02ARMA}.
	
	We remark that the critical and sub-critical spaces play a crucial role in the above mentioned theory for Navier-Stokes equations. However, it is still plausible to study the problem in a large subset of some super-critical space, and one approach among them is the randomized data problem. This subject was initiated by Bourgain \cite{Bourgain-1996, Bourgain-1996-2} for the nonlinear Schr\"odinger equations. Then, Burq and Tzvetkov \cite{Burq-Tzvetkov-2008, Burq-Tzvetkov-2008-2} considered more general random data for the nonlinear wave equations on the compact manifold.
	
	This method has also been applied to Navier-Stokes equations. Zhang-Fang \cite{Zhang-Fang-2011} and  Deng-Cui \cite{Deng-Cui-2011-1} proved the almost sure local well-posedness and  small data global well-posedness  in $L^2(\T^3)$. Later on, Zhang and Fang \cite{Zhang-Fang-2012} further extended the results to the higher dimensional torus $\T^d$ and Euclidean space $\R^d$ when $d \goe3$ with rougher data in $H^{s}$, $-1<s\loe 0$. The same results in $\T^d$ with $d=2$ or $3$ were also studied by Deng and Cui \cite{Deng-Cui-2011-2}. While the above probabilistic results concern the solution for the corresponding integral equation, as for the related weak solution theory, Nahmod, Pavlovi\'{c}, and Staffilani \cite{Nahmod-Pavlovic-Staffilani-2013} showed the almost sure global existence of weak solutions in $H^\al(\mathbb T^d)$ with some $\al<0$ when $d=2, 3$, and they also proved the uniqueness in 2D case. Moreover, the stochastic Navier-Stokes equations were also studied, see for instance \cite{Bensoussan-Temam-1973, Prato-Debussche-2002, Mikulevicius-Rozovskii-2005, Hairer-Mattingly-2006}.
	
	Particularly, in \cite{Zhang-Fang-2012}, the authors introduced a randomization based on a unit-scale decomposition in the frequency space for initial data on $\R^d$. B\'{e}nyi, Oh and Pocovnicu \cite{Benyi-Oh-Pocovnicu-2015} also used similar idea to introduce the Wiener randomization, and the corresponding random data problems for nonlinear Schr\"odinger and wave equations were extensively studied. Furthermore, other kinds of randomization were also introduced in the Euclidean space setting $\mathbb R^d$, such as the ones based on the wave packet \cite{Bringmann-2021}, annuli \cite{Bringmann-2020}, or angular variable \cite{Burq-Krieger-2021} decompositions. Very recently, Shen, Soffer and Wu \cite{Shen-Soffer-Wu-2021-1} constructed a ``narrowed" Wiener randomization, and gave the first probabilistic result  without any regularity restrictions, in the context of nonlinear Schr\"{o}dinger equations.
	
	Furthermore, there are some initial data with very low regularity in practical use, such as the data under the white noise measure in $\dot H^{-d/2-}$ and the Gibbs measure in $\dot H^{-d/2+1-}$. This motivates us to study the problem with arbitrarily rough data, as what was done in \cite{Shen-Soffer-Wu-2021-1}. Therefore, in this note, we aim to apply the ``narrowed" Wiener randomization to study the almost sure well-posedness for the incompressible Navier-Stokes equations. 
	
	\subsection{Definition of randomization}
	
	Before giving the definition, we make a decomposition of the frequency space. Let $a \in \N$. First, for $N \in 2^{\N}$, we denote the cube sets
	$$
	\mathcal{O}_N=\{\xi \in \R^d: |\xi_j|\loe N, j=1,2,\cdots,d\},
	$$
	and then define
	$$
	Q_N=\mathcal{O}_{2N}\backslash\mathcal{O}_N.
	$$
	
	Next, we make a further decomposition of $Q_N$. Note that $a$ is a positive integer, we make a partition of $Q_N$ with the essentially disjoint sub-cubes
	$$
	\mathcal{A}(Q_N):=\{Q: Q  \mbox{ is dyadic cube with length } N^{-a} \mbox{ and }  Q \subset {Q_N}\}.
	$$
	Then, we have $\#\mathcal{A}(Q_N)=(2^d-1)N^{d(a+1)}$ and $Q_N=\cup_{Q \in \mathcal{A}(Q_N)}Q$.
	
	Now, we define the final decomposition
	$$
	\mathcal{Q}:=\{\mathcal{O}_1\} \cup \{Q:Q\in \mathcal{A}(Q_N) \mbox{ and } N\in 2^{\N}\}.
	$$
	By the above construction,
	$$
	\R^d=\mathcal{O}_1 \cup (\cup_{N\in 2^{\N}}Q_N)=\mathcal{O}_1 \cup (\cup_{Q\in \mathcal{A}(Q_N)}Q)=\cup_{Q \in \mathcal{Q}}Q.
	$$
	Therefore, $\mathcal{Q}$ is a countable family of essentially disjoint caps covering $\R^d$. We can renumber the cubes in $\mathcal{Q}$ as follows:
	$$
	\mathcal{Q}=\{Q_j:j \in \N\}.
	$$
	\begin{definition} [``Narrowed'' Wiener Randomization]\label{def-randomization}
		Let $d\goe 2$ and $0<\ep\ll1$. Given $s\in \R$, we set the parameter $a \in \N$ such that $a\goe  \frac{2(d\ep-1-s)}{d(1-2\ep)}$. Let $\widetilde{\psi}_j \in C_0^{\infty}(\R^d)$ be a real-valued function such that $0\loe \widetilde{\psi}_j\loe 1$ and
		$$
		\widetilde{\psi}_j(\xi)=\left\{
		\begin{array}{crl}
			1,    &  &{ when\quad  \xi \in Q_j},\\
			smooth,& &{otherwise},\\
			0,    &  &{ when \quad \xi  \notin 2Q_j},
		\end{array}\right.
		$$
	where we use the notation that $2Q_j$ is the cube with the same center as $Q_j$ and with $\rm{diam}(2Q_j)=2\rm{diam}(Q_j)$. Now, let
	$$
	\psi_j(\xi):=\frac{\widetilde{\psi}_j(\xi)}{\sum_{j' \in \N}\widetilde{\psi}_{j'}(\xi)}.
	$$
	Then, $\psi_j \in  C_0^{\infty}([0, +\infty))$ is a real-valued function, satisfying $0\loe \psi_j\loe 1$ and for all $\xi \in \R^d$, $\sum_{j \in \N}\psi_j (\xi)=1$.

	Denote the Fourier transform on $\R^d$ by $\mathcal{F}$. Then, for any function $f\in H_x^s(\R^d;\R^d)$, we define
	$$
	\Box_jf=\mathcal{F}^{-1}\big(\psi_j (\xi)\mathcal{F}f(\xi)\big).
	$$
	Let $(\Omega, \mathcal{A}, \mathbb{P})$ be a probability space. Let $\{g_{j}(\omega):j\in\N\}$ be a sequence of zero-mean, complex-valued Gaussian random variables on $\Omega$, where the real and imaginary parts of $g_j$ are independent. Then, for the above $f$, we define its randomization $f^{\omega}$ by
	\begin{align*}
		f^{\omega}=\sum_{j\in \N} g_{j}(\omega) \Box_j f.
	\end{align*}
\end{definition}

\subsection{Main results}
In the following, we use the statement ``almost every $\omega \in \Omega$, $PC(\omega)$ holds'' to mean that
$$
\mathbb{P}\big(\{\omega \in \Omega: PC(\omega) \mbox{ holds}\}\big)=1.
$$

In this note, we consider the incompressible Navier-Stokes equations with randomized data:
\begin{equation}\label{NS-randomized data}
	\left\{ \aligned
	&\partial_t u-\Delta u +(u\cdot \nabla) u+\nabla p=0,
	\\
	&\mbox{div } u=0,
	\\
	&u(0,x)=f^{\omega}(x).
	\endaligned
	\right.
\end{equation}
We define the space of solenoidal vector fields $H_\sigma^{s}(\mathbb R^d;\R^d)$ as 
\EQ{
	\overline{\{u \in C_0^\I(\R^d;\R^d): \mbox{div } u=0\}}^{H^{s}}.
}
Our main results are as follows:
\begin{thm}\label{main theorem}
	Let $d\goe2$ and $0<\ep\ll1$. Given any $s\in \R$ and $f \in H_\sigma^{s}(\R^d;\R^d)$. Suppose that the randomization $f^{\omega}$ is defined in Definition \ref{def-randomization}. Then, for almost every $\omega \in \Omega$, we have that for any $2<q<+\I$,
	\EQ{
	f^\om\in \dot B^{d\ep-1}_{1/\ep,q}(\R^d;\R^d).
	}
	Hence, the following results hold:
	\vspace{-0.2cm}
	\begin{enumerate}
		\item 
		There exists $T=T(f^{\omega})>0$, such that \eqref{NS-randomized data} is locally well-posed on $[0,T]$.
		\item
		If $d=2$, then  \eqref{NS-randomized data} is globally well-posed.
	\end{enumerate}
\end{thm}

We make some remarks regarding our results. Even though this randomization cannot give the smoothing effect in the $L_x^2$-based Sobolev's space, it indeed implies that the initial data belongs to some critical Besov space. Note that in such Besov space, Cannone \cite{Cannone-1997} proved the local well-poseness for $d\goe2$, and Gallagher and Planchon \cite{GP02ARMA} proved 2D global result, then the related well-posedness holds immediately. Since the proof is too trivial in our setting, we will not submit this note for publication.

In our previous version, we prove the local well-posedness by improving the nonlinear part $u-e^{t\De}f^\om$ in higher order Sobolev space $\dot H^{\frac d2-1}(\R^d)$. We think that this improvement is not so interesting in the context of Navier-Stokes equations, because the well-posedness theory in $L^p$-based spaces is well understood, as described above.

Our result shows that the condition
\EQ{
a\goe \frac{2(d\ep-1-s)}{d(1-2\ep)}=-\frac2d(1+s)+
}
guarantees the well-posedness of Navier-Stokes equations. Particularly, if $-1<s<0$, we can take $a=0$, which corresponds to the classical Wiener type randomization.

Moreover, if $a$ is too large, then the initial data can be smoother in the $C^k$ sense. For example, if we take $a > \frac2d(k-s)$ for some large $k\in\N^+$, we can similarly prove that
\EQ{
\norm{\jb{\nabla}^{k}\Box_j f}_{L_x^{\I-}(\R^d)} \lsm \norm{\Box_jf}_{H_x^s(\R^d)},
}
which implies $f^\om\in W^{k,\I-}$ almost surely. Then by Morrey's inequality, we obtain $f^\om\in C^{k-1,\ga}$ for some $0<\ga<1$ and a.e. $\om\in\Om$.

\section{Preliminary}

\subsection{Notations}
We denote $\R^+=[0,+\I)$. For any $a\in \R$, $a\pm :=a \pm \epsilon$ for arbitrary small $\epsilon >0$. For any set $A$, we denote $\# A$ as the cardinal number of $A$. We write $X\lesssim Y$ or $Y\gtrsim X$ to denote the estimate $X\loe CY$ for some constant $C>0$. Throughout the whole paper, the letter $C$ will denote different positive constants which are not important in our analysis and may vary line by line. If $f\loe Cg$ and $g\loe Cf$, we write $f\thicksim g$.
Moreover, we write ``$\mbox{a.e. }\omega \in \Omega$" to mean that ``almost every $\omega \in \Omega$".


Now, we give the definitions of various vector-valued functional spaces. Given $1\loe p \loe \infty$, $s\in\R$, and $d, n\in\N^+$. Let $\mathcal{S}(\R^d;\R^n)$ be the Schwartz space, $\mathcal{S'}(\R^d;\R^n)$ be the tempered distribution space, and $C_0^\I(\mathbb R^d;\R^n)$ be the space of all the smooth compact-supported functions. $L^p(\R^d;\R^n)$ denotes the usual Lebesgue space.  $W^{s,p}(\R^d;\R^n)$ and $\dot W^{s,p}(\R^d;\R^n)$ denote the inhomogeneous and homogeneous Sobolev spaces, respectively. We denote the related Bosev spaces as $B^s_{p,q}(\R^d;\R^n)$ and $\dot B^s_{p,q}(\R^d;\R^n)$.
We also denote that $H^s(\R^d;\R^n):=W^{s, 2}(\R^d;\R^n)$ and $\dot{H}^s(\R^d;\R^n):=\dot{W}^{s, 2}(\R^d;\R^n)$.
For the functional spaces $X$ defined above, we define
\EQ{
	X_\sigma(\R^d;\R^d) := \overline{\{u \in C_0^\I(\mathbb R^d;\R^d):\mbox{div } u=0\}}^{X}.
}
Moreover, if $d=n$, we denote that $X(\R^d):=X(\R^d;\R^d)$. We also write $X(\R^d)$ as $X$ for short. Finally, we use the following norms to denote the mixed spaces $L_t^qL_x^r$, that is
\begin{align*}
	\|u\|_{L_t^qL_x^r}=\big(\int \|u(t,\cdot)\|_{L_x^r}^qdt\big)^{\frac{1}{q}}.
\end{align*}
When $q=r$, we abbreviate $L_t^qL_x^q$ to $L_{t,x}^q$.

\subsection{Basic lemmas}
In this subsection, we gather some useful results. The followings are referred to Lemma 2.6 and Lemma 2.7 in \cite{Shen-Soffer-Wu-2021-1}.

\begin{lem}[Orthogonality]\label{Orthogonality}
	Let $f\in L_x^2(\R^d;\R)$. Then, we have
	\begin{align*}
		\|\Box_jf\|_{L_x^2l_{j\in \N}^2}\thicksim \|f\|_{L_x^2}.
	\end{align*}
\end{lem}

\begin{lem}[$L^q$-$L^p$ estimate]\label{lem:bernstein}
	Let $a>0$ and $2\loe p\loe q \loe\infty$. Given any $j\in\N$, then
	\EQ{
		\norm{\Box_j f}_{L_x^q(\R^d)}\lesssim \norm{\jb{\nabla}^{-a(\frac dp-\frac dq)}\Box_j  f}_{L_x^p(\R^d)}.
	}
\end{lem}

We also recall some probabilistic results.
\begin{lem}[Large deviation estimate, \cite{Burq-Tzvetkov-2008}]\label{Large derivation}
	
	Let $(\Omega, \mathcal{A}, \mathbb{P})$ be a probability space. Let $\{g_j\}_{j \in \mathbb{N}^+}$ be a sequence of real-valued, independent, zero-mean random variables with associated distributions $\{\mu_j\}_{j \in \mathbb{N}^+}$ on $\Omega$. Suppose $\{\mu_j\}_{j \in \mathbb{N}^+}$ satisfies that there exists $c>0$ such that for all $\gamma \in \R$ and $j \in \mathbb{N}^+$
	\begin{align*}
		\big|\int_{\R} e^{\gamma x}d\mu_j(x)\big|\loe e^{c\gamma^2},
	\end{align*}
	then there exists $\alpha >0$ such that for any $\lambda>0$ and any complex-valued sequence $\{c_j\}_{j \in \mathbb{N}^+} \in l_j^2$, we have
	\begin{align*}
		\mathbb{P}(\{\omega:\big|\sum_{j=1}^{\infty}c_j g_j(\omega)\big|>\lambda\})\loe 2 e^{-\alpha \lambda \|c_j\|_{l_j^2}^{-2}}.
	\end{align*}
	Furthermore, there exists $C>0$ such that for any $2\loe p < \infty$ and complex-valued sequence $\{c_j\}_{j \in \mathbb{N}^+} \in l_j^2$, we have
	\begin{align*}
		\big\|\sum_{j=1}^{\infty}c_j g_j(\omega)\big\|_{L_\omega^p(\Omega)}\loe C \sqrt{p}\|c_j\|_{l_j^2}.
	\end{align*}
\end{lem}

The following lemma can be proved by the method in \cite{Tzvetkov-2010}.

\begin{lem}\label{large P thing}
	Let $F$ be a real-valued measurable function on a probability space $(\Omega, \mathcal{A}, \mathbb{P})$. Suppose that there exist $C_0>0$, $K>0$ and $p_0\goe 1$ such that for any $p\goe p_0$, we have
	\begin{align*}
		\|F\|_{L_\omega^p(\Omega)}\loe\sqrt{p}C_0K.
	\end{align*}
	Then, there exist $c>0$ and $C_1>0$, depending on $C_0$ and $p_0$ but independent of $K$, such that for any $\lambda>0$,
	\begin{align*}
		\mathbb{P}(\{\omega\in \Omega:|F(\omega)|>\lambda\})\loe C_1 e^{-c \lambda^2 K^{-2}}.
	\end{align*}
\end{lem}

\section{A short proof of the main results}

Let $p=1/\ep$ and $2<q<+\I$. Now, we restrict the variables on $x\in\R^d$, $\om\in\Om$, $j\in\N$, and $N\in2^\Z$. Define $P_N$ as the dyadic Littlewood-Paley projection operator. For any $N\goe 1$, since $a\goe \frac{2(d\ep-1-s)}{d(1-2\ep)}$, by Lemma \ref{lem:bernstein},  we have for any $j\in \N$, 
\EQ{
\normb{|\nabla|^{\frac dp-1}P_N\Box_j f}_{L_x^p} \lsm \normb{|\nabla|^{\frac dp-1}\jb{\nabla}^{-a(\frac d2-d\ep)}P_N\Box_j f}_{L_x^2} \lsm \norm{P_N\Box_jf}_{H_x^s}.
}
For any $N\loe 2^{-1}$, in fact $j=0$ in this case, then by the usual Bernstein inequality for the operator $P_N$,
\EQ{
\normb{|\nabla|^{\frac dp-1}P_N\Box_j f}_{L_x^p} \lsm N^{\frac dp-1+\frac d2-\frac dp}\normb{P_N\Box_0 f}_{L_x^2} \lsm \norm{P_N\Box_0f}_{H_x^s}.
}
Therefore, for any $N\in2^\Z$ and $j\in\N$,
\EQn{\label{eq:bernstein}
\normb{|\nabla|^{\frac dp-1}P_N\Box_j f}_{L_x^p} \lsm \norm{P_N\Box_jf}_{H_x^s}.
}
Then, 
For any $\rho\goe p_0= \max\fbrk{p,q}$, by Minkowski's inequality, Lemma \ref{Large derivation}, \eqref{eq:bernstein}, Lemma \ref{Orthogonality}, and $l_N^2\hra l_N^q$,
\EQn{\label{eq:initial-data-1}
	\norm{f^\om}_{L_\om^\rho \dot B^{\frac dp-1}_{p,q}}\lesssim & \normb{|\nabla|^{\frac dp-1}P_Nf^\om}_{l_N^qL_x^pL_\om^\rho} \\ \lsm & \sqrt{\rho} \normb{|\nabla|^{\frac dp-1}P_N\Box_jf}_{l_N^qL_x^pl_j^2} \\
	\lsm & \sqrt{\rho} \normb{|\nabla|^{\frac dp-1}P_N\Box_jf}_{l_j^2l_N^qL_x^p} \\
	\lsm & \sqrt{\rho} \norm{P_N\Box_jf}_{l_j^2l_N^2H_x^s} \\
	\lsm & \sqrt{\rho} \norm{f}_{H_x^s}.
}
For any $2\loe\rho< p_0= \max\fbrk{p,q}$, by H\"older's inequality in $\om$,
\EQ{
	\norm{f^\om}_{L_\om^\rho \dot B^{\frac dp-1}_{p,q}} \lsm \norm{f^\om}_{L_\om^{p_0} \dot B^{\frac dp-1}_{p,q}}.
}
Then, using \eqref{eq:initial-data-1},
\EQ{
	\norm{f^\om}_{L_\om^{p_0} \dot B^{\frac dp-1}_{p,q}} \lsm \sqrt{p_0} \norm{f}_{H_x^s}\lsm_{p,q} \sqrt{\rho} \norm{f}_{H_x^s}.
}
This gives for any $2\loe \rho<p_0$,
\EQn{\label{eq:initial-data-2}
	\norm{f^\om}_{L_\om^\rho \dot B^{\frac dp-1}_{p,q}} \lsm \sqrt{\rho} \norm{f}_{H_x^s}.
}
Now, combining \eqref{eq:initial-data-1} and \eqref{eq:initial-data-2}, for any $\rho\goe 2$,
\EQ{
\norm{f^\om}_{L_\om^\rho \dot B^{\frac dp-1}_{p,q}} \lsm & \sqrt{\rho} \norm{f}_{H_x^s},
}
then by Lemma \ref{large P thing}, we obtain for almost every $\om\in\Om$,
\EQ{
\norm{f^\om}_{\dot B^{\frac dp-1}_{p,q}}<+\I.
}
Then, the local well-posedness for $d\goe 2$ follows by Cannone's result \cite{Cannone-1997}, and 2D global well-posedness follows by Gallagher and Planchon's result \cite{GP02ARMA}.

\bigskip
\section*{Acknowledgment}
J. Shen is partially supported by NSFC 12171356. The authors are very grateful to the anonymous referee for valuable comments, and for pointing out that our results can be proved easily. The authors would also like to thank Prof. Yifei Wu and Prof. Kenji Nakanishi for helpful discussions and valuable suggestions.


\begin{thebibliography}{99}
	
	
	
	
	
	\bibitem{Bensoussan-Temam-1973}
	A. Bensoussan and R. Temam, \`{E}quations stochastiques du type Navier-Stokes, J. Funct. Anal., 13 (2), 1973, 195-222.
	
	
	
	
	\bibitem{Benyi-Oh-Pocovnicu-2015}
	A. B\'{e}nyi, T. Oh and O. Pocovnicu,  On the probabilistic Cauchy theory of the cubic nonlinear Schr\"{o}dinger equation on $\R^d$, $d\geq3$, Trans. Amer. Math. Soc.
	Ser. B, 2 (2015), 1-50.
	
	
	\bibitem{Bourgain-1996}
	J. Bourgain, Periodic nonlinear Schr\"{o}dinger equation and invariant measures,
	Comm. Math. Phys., 166 (1), 1994, 1-26.
	
	
	
	\bibitem{Bourgain-1996-2}
	J. Bourgain, Invariant measures for the 2D-defocusing nonlinear Schr\"{o}dinger
	equation, Comm. Math. Phys., 176 (2), 1996, 421-445.
	
	
	\bibitem{Bourgain-Pavlovic-2008}
	J. Bourgain and N. Pavlovi\'{c}, Ill-posedness of the Navier-Stokes equations in a critical space in 3D, J. Funct. Anal., 255 (9), 2008, 2233-2247.
	
	
	
	\bibitem{Bringmann-2020}
	B. Bringmann, Almost-sure scattering for the radial energy-critical nonlinear wave equation in three dimensions, Anal. PDE, 13 (4), 2020, 1011-1050.
	
	
	\bibitem{Bringmann-2021}
	B. Bringmann, Almost sure scattering for the energy critical nonlinear wave equation, Amer. J. Math., 143 (6), 2021, 1931-1982.
	
	\bibitem{Buckmaster-Vicol-2019}
	T. Buckmaster and V. Vicol, Nonuniqueness of weak solutions to the Navier-Stokes equation, Ann. of Math. (2), 189 (1), 2019, 101-144.
	
	
	\bibitem{Burq-Krieger-2021}
	N. Burq and J. Krieger, Randomization improved Strichartz estimates and global well-posedness for supercritical data, Ann. Inst. Fourier, 2021. DOI: 10.5802/aif.3448.
	
	\bibitem{Burq-Tzvetkov-2008}
	N. Burq and N. Tzvetkov, Random data Cauchy theory for supercritical wave equations. I. Local theory, Invent. Math., 173 (3), 2008, 449-475.
	
	
	\bibitem{Burq-Tzvetkov-2008-2}
	N. Burq and N. Tzvetkov, Random data Cauchy theory for supercritical wave equations. II. A global existence result, Invent. Math., 173 (3), 2008,
	477-496.
	
	
	\bibitem{Caffarelli-Kohn-Nirenberg-1982}
	L. Caffarelli, R. Kohn and L. Nirenberg, Partial regularity of suitable weak solutions of the Navier-Stokes equations, Communications on pure and applied mathematics, 35 (6), 1982, 771-831.
	
	\bibitem{Cannone-1997}
	M. Cannone, A generalization of a theorem by Kato on Navier-Stokes equations, Rev. Mat. Iberoamericana, 13 (3), 1997, 515-541.
	
	\bibitem{Cannone-Planchon-1996}
	M. Cannone and F. Planchon, Self-similar solutions for Navier-Stokes equations in $\R^3$, Comm. Partial Differential Equations, 21 (1-2), 1996, 179-193.
	
	
	
	
	\bibitem{Chemin-1992}
	J. Y. Chemin, Remarques sur l'existence globale pour le syst\`{e}me de Navier-Stokes incompressible, SIAM J. Math. Anal., 23 (1), 1992, 20-28.
	
	\bibitem{Chemin-1999}
	J. Y. Chemin, Th\'eor\`{e}mes d'unicit\'e pour le syst\`eme de Navier-Stokes tridimensionnel. (French) [Uniqueness theorems for the three-dimensional Navier-Stokes system],
	J. Anal. Math., 77 (1999), 27-50.
	
	
	\bibitem{Chemin-Gallagher-2009}
	J. Y. Chemin and I. Gallagher, Wellposedness and stability results for the Navier-Stokes equations in $\R^3$, Ann. Inst. H. Poincare Anal. Non Lineaire, 26 (2), 2009, 599-624.
	
	
	
	
	
	\bibitem{Deng-Cui-2011-1}
	C. Deng and S. Cui, Random-data Cauchy problem for the Navier-Stokes equations on $\T^3$ , J. Differential Equations, 251 (4-5), 2011, 902-917.
	
	
	\bibitem{Deng-Cui-2011-2}
	C. Deng and S. Cui, Random-data Cauchy problem for the periodic Navier-Stokes equations with initial data in negative-order Sobolev spaces, arXiv:1103.6170v1.
	
	
	
	
	
	\bibitem{Escauriaza-Seregin-Sverak-2003}
	L. Escauriaza, G. A. Seregin and V. \v{S}ver\'{a}k, $L_{3,\I}$-solutions of Navier-Stokes equations and backward uniqueness, Uspekhi Mat. Nauk 58, 350 (2), 2003, 3–44.
	
	
	\bibitem{Fujita-Kato-1964}
	H. Fujita and T. Kato, On the Navier-Stokes initial value problem. I, Arch. Rational Mech. Anal., 16 (1964), 269-315.
	
	
	
	\bibitem{GP02ARMA}
	I. Gallagher and F. Planchon,
	\newblock{On global infinite energy solutions to the Navier-Stokes equations in two dimensions,} 
	Arch. Ration. Mech. Anal., 161 (4), 2002, 307--337.
	
	
	\bibitem{Hairer-Mattingly-2006}
	M. Hairer and J. C. Mattingly, Ergodicity of the 2D Navier-Stokes equations with degenerate stochastic forcing, Ann. of Math. (2), 164 (3), 2006, 993-1032.
	
	
	\bibitem{Hopf-1951}
	E. Hopf, \"{U}ber die Anfangswertaufgabe f\"{u}r die hydrodynamischen Grundgleichungen, (German) Math. Nachr., 4 (1951), 213-231.
	
	
	
	\bibitem{Jia-Sverak-2014}
	H. Jia and V. \v{S}ver\'{a}k, Local-in-space estimates near initial time for weak solutions of the Navier-Stokes equations and forward self-similar solutions, Invent. Math., 196 (1), 2014, 233-265.
	
	
	
	\bibitem{Kato-1984}
	T. Kato, Strong $L^p$ solutions of the Navier-Stokes equation in $\R^m$, with applications to weak solutions, Math. Z., 187 (4), 1984, 471-480.
	
	
	
	
	
	
	
	\bibitem{Koch-Tataru-2001}
	H. Koch and D. Tataru, Well-posedness for the Navier-Stokes equations, Adv. Math., 157 (2001), 22-35.
	
	
	
	\bibitem{Koch-Nadirashvili-Seregin-Sverak-2009}
	G. Koch, N. Nadirashvili, G. A. Seregin and  V. \v{S}ver\'{a}k, Liouville theorems for the Navier-Stokes equations and applications, Acta Math., 203 (1), 2009, 83-105.
	
	
	
	\bibitem{Ladyzhenskaya-1963}
	O. A. Ladyzhenskaya, The mathematical theory of viscous incompressible flow, Gordon and Breach, New York-London, 1963.
	
	
	\bibitem{Leray-1934}
	J. Leray, Sur le mouvement d'un liquide visqueux emplissant l'espace, (French) Acta Math, 63 (1), 1934, 193-248.
	
	\bibitem{Lei-Lin-2011}
	Z. Lei and F. H. Lin, Global mild solutions of Navier-Stokes equations, Comm. Pure Appl. Math., 64 (2011), 1297-1304.
	
	\bibitem{Lemarie-Rieusset-2002}
	P. G. Lemarie-Rieusset, Recent Developments in the Navier-Stokes Problem, Chapman and Hall/CRC Research Notes in Mathematics,
	vol. 431, Chapman and Hall/CRC, Boca Raton, FL, 2002, 395 p.
	
	
	
	
	
	\bibitem{Lin-1998}
	F. H. Lin, A new proof of the Caffarelli-Kohn-Nirenberg theorem, Communications on Pure and Applied Mathematics, 51 (3), 1998, 241-257.
	
	
	
	\bibitem{Mikulevicius-Rozovskii-2005}
	R. Mikulevicius and B. L. Rozovskii, Global $L_2$-solutions of stochastic Navier-Stokes equations, Ann. Probab., 33 (1), 2005, 137-176.
	
	
	\bibitem{Nahmod-Pavlovic-Staffilani-2013}
	A. Nahmod, N. Pavlovi\'{c} and G. Staffilani, Almost sure existence of global weak solutions for super-critical Navier-Stokes equations, SIAM J. Math. Anal., 45 (6), 2013, 3431-3452.
	
	\bibitem{NRS96Acta}
	J. Nečas,  M. Růžička and V. Šverák, On Leray's self-similar solutions of the Navier-Stokes equations, Acta Math., 176 (2), 1996, 283–294.
	
	
	\bibitem{Prato-Debussche-2002}
	G. Da Prato and A. Debussche, Two-dimensional Navier-Stokes equations driven by a space-time white noise, J. Funct. Anal., 196 (1), 2002, 180-210.
	
	
	
	\bibitem{Scheffer-1977}
	V. Scheffer, Hausdorff measure and the Navier-Stokes equations, Communications in Mathematical Physics, 55 (2), 1977, 97-112.
	
	
	\bibitem{Shen-Soffer-Wu-2021-1}
	J. Shen, A. Soffer and Y. Wu, Almost sure scattering for the nonradial energy-critical NLS with arbitrary regularity in 3D and 4D cases, arXiv:2111.11935.
	
	\bibitem{Tsa98ARMA}
	T. Tsai, On Leray's self-similar solutions of the Navier-Stokes equations satisfying local energy estimates, Arch. Rational Mech. Anal., 143 (1), 1998, 29–51.
	
	\bibitem{Tzvetkov-2010}
	N. Tzvetkov, Construction of a Gibbs measure associated to the periodic Benjamin-Ono equation, Probab. Theory Related Fields, 146 (3-4), 2010, 481-514.
	
	\bibitem{Wang15Adv}
	B. Wang, Ill-posedness for the Navier-Stokes equations in critical Besov spaces $\dot B^{-1}_{\I,q}$, Adv. Math., 268 (2015), 350–372.
	
	
	\bibitem{Zhang-Fang-2011}
	T. Zhang and D. Fang, Random data Cauchy theory for the incompressible three dimensional Navier-Stokes equations, Proc. AMS., 139 (8), 2011, 2827-2837.
	
	
	\bibitem{Zhang-Fang-2012}
	T. Zhang and D. Fang, Random data Cauchy theory for the generalized incompressible Navier-Stokes equations, J. Math. Fluid Mech., 14 (2), 2012, 311-324.
	
	
	
	
	
	
	
	
	
	
	
	
	
\end{thebibliography}
\end{document}